\documentclass[12pt]{amsart}
\usepackage{graphicx}
\newtheorem{theorem}{Theorem}[section]

\newtheorem{corollary}[theorem]{Corollary}

\theoremstyle{remark}
\newtheorem{remark}[theorem]{Remark}

\numberwithin{equation}{section}

\begin{document}

\title[
Irregular smoothing and the number of Reidemeister moves
]{
Irredular smoothing and the number of Reidemeister moves
}

\author{Chuichiro Hayashi and Miwa Hayashi}


\thanks{The first author is partially supported
by Grant-in-Aid for Scientific Research (No. 22540101),
Ministry of Education, Science, Sports and Technology, Japan.}

\begin{abstract}
 In the previous paper \cite{HH2},
we consider a link diagram invariant of Hass and Nowik type
using regular smoothing and unknotting number,
to estimate the number of Reidemeister moves needed for unlinking.
 In this paper,
we introduce a new link diagram invariant
using irregular smoothing,
and give an example of a knot diagram of the unknot
for which the new invariant gives a better estimation
than the old one.
\end{abstract}

\maketitle

\section{Introduction}\label{section:introduction}

\begin{figure}[htbp]
\begin{center}
\includegraphics[width=9cm]{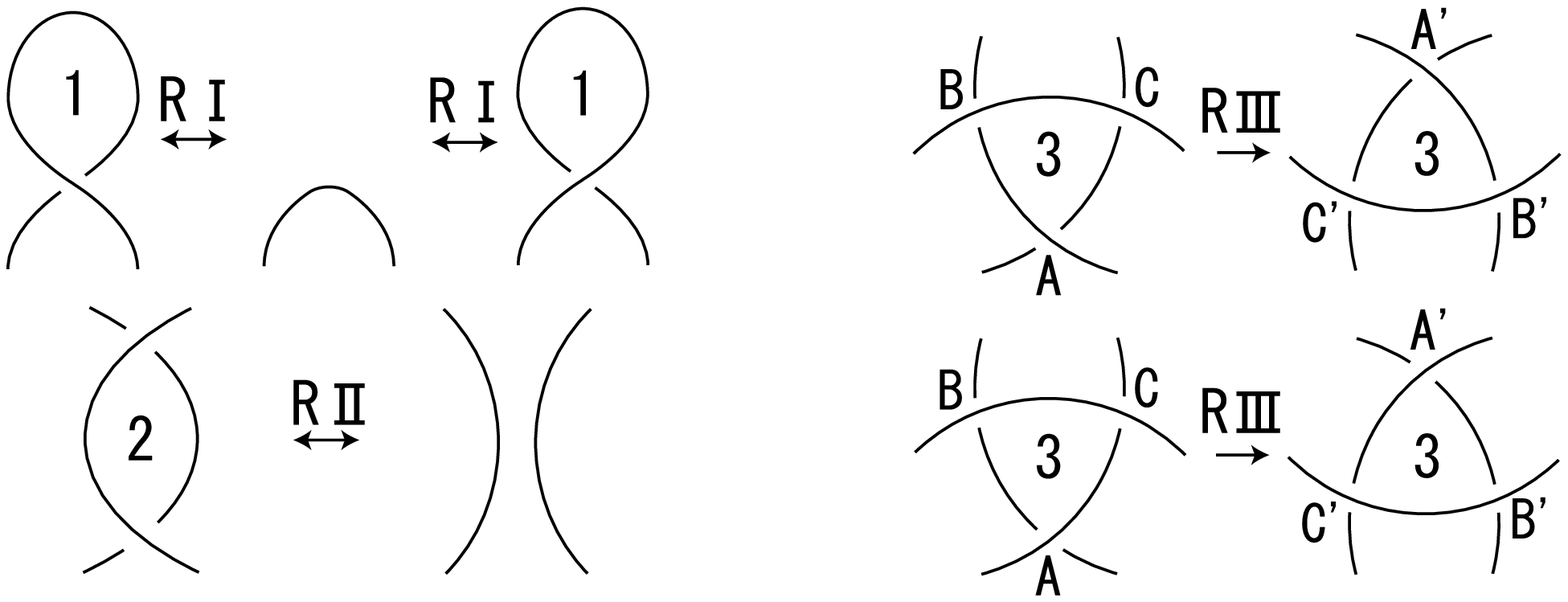}
\end{center}
\caption{}
\label{fig:Reid123}
\end{figure}

\begin{figure}[htbp]
\begin{center}
\includegraphics[width=5cm]{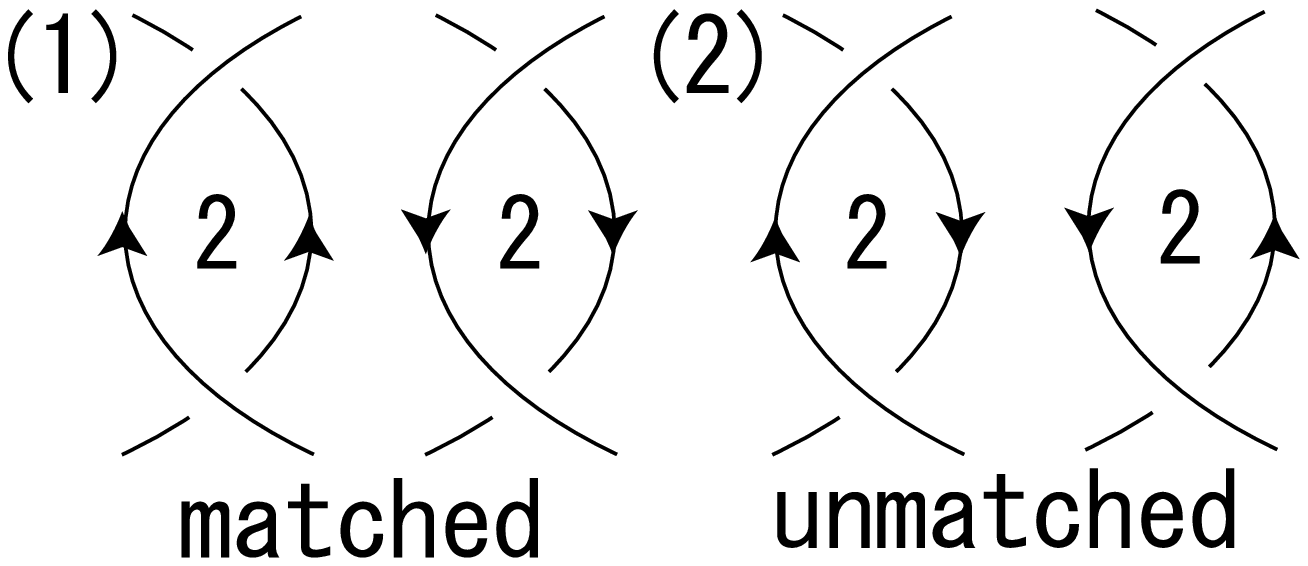}
\end{center}
\caption{}
\label{fig:matched}
\end{figure}

 In this paper,
we regard a knot as a link with one component,
and assume that link diagrams are in the $2$-sphere. 
 Links and link diagrams are assumed to be oriented
unless otherwise specified.
 A {\it Reidemeister move} is a local move of a link diagram
as in Figure \ref{fig:Reid123}.
 An RI (resp. II) move
creates or deletes a monogon face (resp. a bigon face).
 An RII move is called {\it matched} and {\it unmatched} 
according to the orientations of the edges of the bigon
as shown in Figure \ref{fig:matched}.
 An RIII move is performed on a $3$-gon face,
deleting it and creating a new one.
 Any such move does not change the link type.
 As Alexander and Briggs \cite{AB} and Reidemeister \cite{R} showed,
for any pair of diagrams $D_1$, $D_2$ which represent the same link type,
there is a finite sequence of Reidemeister moves
which deforms $D_1$ to $D_2$.

 Necessity of Reidemeister moves of type II and III is studied 
in \cite{O}, \cite{Ma} and \cite{Hg}.
 There are several studies 
of lower bounds for the number of Reidemeister moves
connecting two knot diagrams of the same knot.
 See \cite{SS}, \cite{H}, \cite{CESS}, \cite{HN1}, \cite{HN2}, 
\cite{HH1}, \cite{HHSY} and \cite{HH2}.
 Hass and Nowik introduced in \cite{HN1}
a certain knot diagram invariant $I_{\phi}$
by using (regular) smoothing and an arbitrary link invariant $\phi$. 
 In \cite{HN2}, setting $\phi$ to be the linking number,
they showed that a certain infinite sequence of diagrams of the unknot
requires quadratic number of Reidemeister moves with respect to the number of crossings 
for being unknotted.
 In the previous paper \cite{HH2}, 
we consider a link diagram invariant of Hass and Nowik type
using the unknotting number.
 Note that we obtain a diagram of a knot rather than a multi-component link
when we perform a smooting at a crossing
between ditinct components of a $2$-component link.
 We have considered the unknotting number rather than the linking number
to deal with this situation.
 This also allows us to perform an irregular smoothing on a knot diagram.
 In this paper,
we introduce a new link diagram invariant using irregular smoothing,
and give an example of a diagram of the unknot
for which the new invariant gives 
a better estimation of the number of Reidemeister moves
than the old one.

\begin{figure}[htbp]
 \begin{minipage} {0.4\hsize}
  \begin{center}
   \includegraphics[width=60mm]{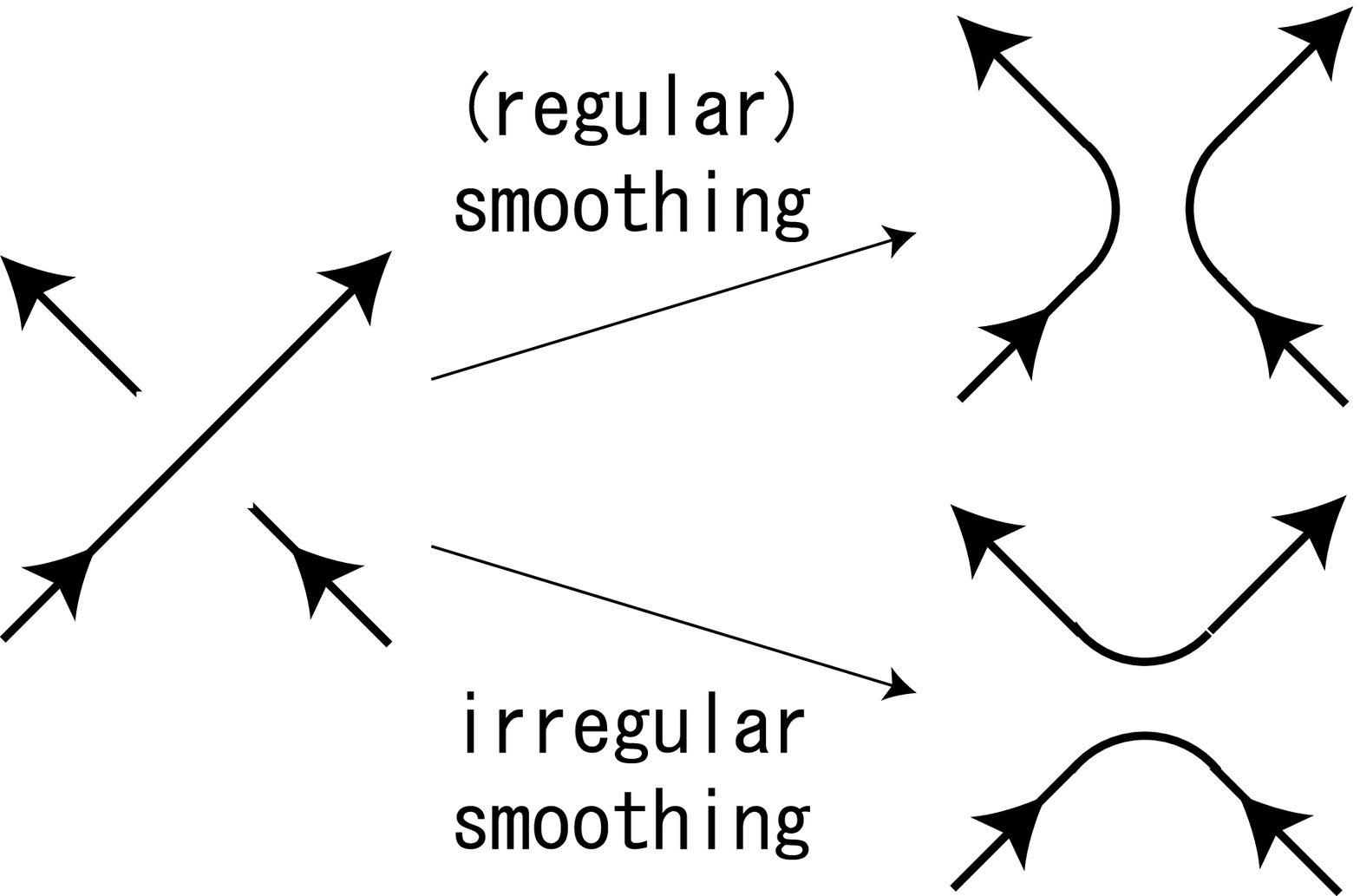}
  \end{center}
  \caption{}
  \label{fig:smoothing}
 \end{minipage}
 \begin{minipage} {0.4\hsize}
  \begin{center}
   \includegraphics[width=60mm]{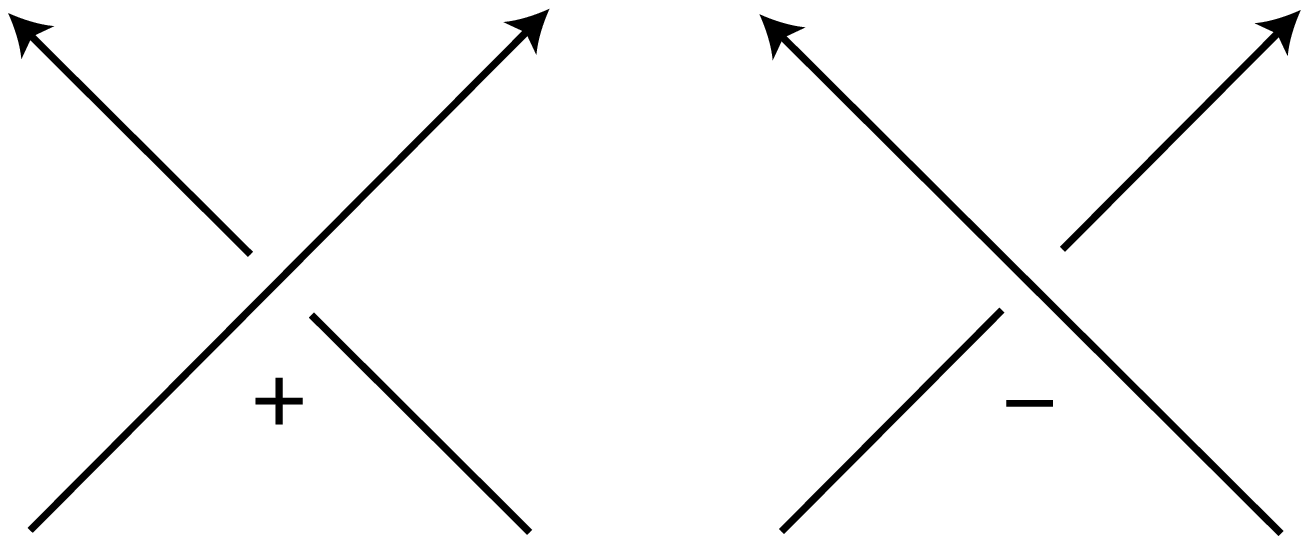}
  \end{center}
  \caption{}
  \label{fig:sign}
 \end{minipage} 
\end{figure}

 An {\it irregular smoothing} at a crossing of a link diagram 
is an operation as in Figure \ref{fig:smoothing}.
 We regard the resulting link digram to be unoriented.
 For a link $K$ of $m$ components,
the {\it unknotting number} $u(K)$ of $K$ is 
the $X$-Gordian distance 
between $K$ and the trivial link of $m$ components.
 (Precise descriptions of the definitions 
of the irregular smoothing and the unknotting number are given 
in Section \ref{section:change}.)
 Then we define the new link diagram invariant $iu_{S,T}$ as below.
 Let $L$ be an oriented $n$-component link
with its components numbered as $L_1, L_2, \cdots L_n$.
 Let $D$ be a diagram of $L$ in the $2$-sphere.
 For a crossing $x$ of $D$,
let $D_x$ (resp. $\check{D_x}$) denote the link obtained from $D$ 
by performing a regular smoothing operation 
(resp. an irregular smoothing operation) at $x$.
 Note that $D_x$ ($\check{D_x}$) is a link rather than a diagram.
 The link $D_x$ (resp. $\check{D_x}$) has $n+1$ (resp. $n$) components 
when $x$ is a crossing between subarcs of the same component,
and $n-1$ (resp. $n-1$) components 
when $x$ is a crossing between subarcs of distinct components.
 Let $S$ and $T$ be symmetric matrices of order $n$ 
with each element $s_{ij}$, $t_{ij}$ ($1 \le i, j \le n$) 
is equal to $-1$, $0$ or $+1$
such that $t_{ij}=0$ if and only if $s_{ij}=0$.
 We set ${\mathcal C}(D)$ (resp. $\check{\mathcal C} (D)$) 
to be the set of all the crossings of $D$
between subarcs of components $L_i$ and $L_j$
with the $(i,j)$-element $s_{ij} = +1$ (resp. $s_{ij}=-1$).
 Then we set 
$$iu_{S,T}(D) = 
\sum_{x \in {\mathcal C}(D)} t(x)\cdot{\rm sign}(x)\cdot|\Delta u(D_x)|
+ \sum_{x \in \check{\mathcal C}(D)} t(x)\cdot{\rm sign}(x)
\cdot|\Delta u(\check{D_x})|,$$
where
$\Delta u(D_x)$ (resp. $\Delta u(\check{D_x})$) is the difference 
between the unknotting numbers of $D_x$ and $L$ (resp. $\check{D_x}$ and $L$),
i.e., $\Delta u (D_x)= u(D_x) - u(L)$ 
(resp. $\Delta u (\check{D_x})= u(\check{D_x}) - u(L)$),
and $t(x)$ denotes the $(i,j)$-element of $T$
such that $x$ is a crossing between a subarc of the component $L_i$ and that of $L_j$.
 The sign of a crossing ${\rm sign}(x)$ 
is defined as in Figure \ref{fig:sign} as usual.
 We set $iu_S(D)=0$ for a diagram $D$ with no crossings.

\begin{figure}[htbp]
\begin{center}
\includegraphics[width=14cm]{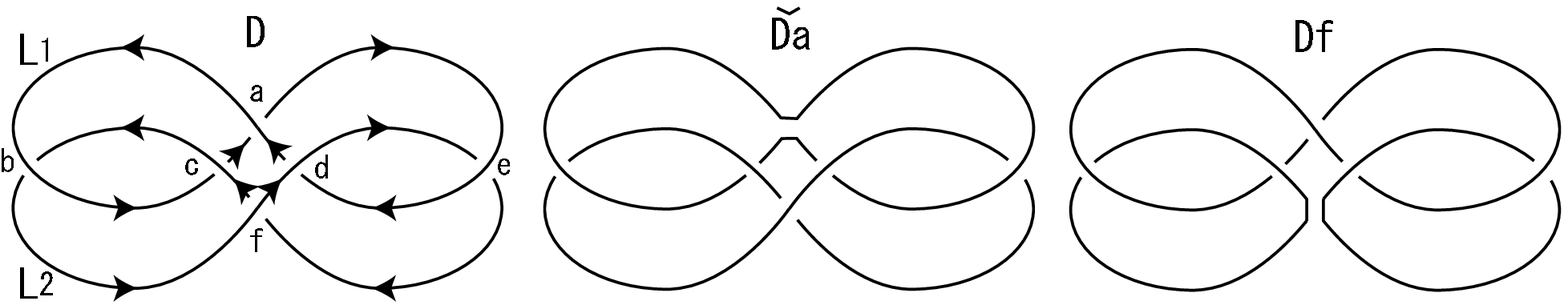}
\end{center}
\caption{}
\label{fig:example}
\end{figure}

 For example,
let $D$ be a link diagram in Figure \ref{fig:example}.
 The link has two components $L_1$ and $L_2$ and is trivial.
 We set 
$S = T = \left( \begin{array}{cc} -1 & 0 \\ 0 & 1 \end{array} \right)$.
 Then the irregular smoothing at the crossing $a$ yields $\check{D_a}$
which is the $(2,4)$-torus link 
and has the unknotting number $u(\check{D_a})=2$.
 Since $t_{11}=-1$ and ${\rm sign}(a)=-1$,
the crossing $a$ contributes by $(-1)(-1)|2-0|$ to $iu_{S,T}(D)$.
 The crossings $b$ through $e$ do not contribute to $iu_{S,T}(D)$
because $s_{12}=0$.
 The link $D_f$ obtained from $D$ by the regular smoothing at $f$ 
is a composite of two Hopf links,
and has the unknotting number $u(D_f)=2$.
 Since $t_{22}=+1$ and ${\rm sign}(f)=+1$,
the crossing $f$ contributes by $(+1)(+1)|2-0|$ to $iu_{S,T}(D)$.
 Hence we have 
$iu_{S,T}(D)= (-1)(-1)|2-0|+(+1)(+1)|2-0|=4$.

\begin{remark}
 If we replace 
$|\Delta u(D_x)|$ and $|\Delta u(\check{D_x})|$
with 
$\Delta u(D_x)$ and $\Delta u(\check{D_x})$,
in the definition of $iu_{S,T}(D)$,
we obtain another link diagram invariant as below.
$$\tilde{iu}_{S,T}(D) = 
\sum_{x \in {\mathcal C}(D)} t(x)\cdot{\rm sign}(x)\cdot\Delta u(D_x)
+ \sum_{x \in \check{\mathcal C}(D)} t(x)\cdot{\rm sign}(x)
\cdot\Delta u(\check{D_x})$$
 All the results about $iu_{S,T}(D)$ in this paper 
hold also for $\tilde{iu}_{S,T}(D)$.
\end{remark}

\begin{theorem}\label{theorem:change}
 The link diagram invariant $iu_{S,T}(D)$ does not change under an RI move,
and changes at most by one under an RII move,
and at most by two under an RIII move.
\end{theorem}

 The above theorem is proved in Section \ref{section:change}.

\begin{corollary}\label{corollary:change1}
 Let $D_1$ and $D_2$ be link diagrams of the same oriented link.
 We need at least $\lceil | iu_{S,T}(D_1) - iu_{S,T}(D_2) | /2 \rceil$ 
RII and RIII moves
to deform $D_1$ to $D_2$ by a sequence of Reidemeister moves,
where $\lceil x \rceil$ denotes the smallest integer 
larger than or equal to $x$.
 In particular,
when $D_2$ is a link diagram 
with ${\mathcal C}(D_2) \cup \check{{\mathcal C}}(D_2) = \emptyset$,
we need at leat $\lceil |iu_{S,T}(D_1)|/2 \rceil$ RII and RIII moves.
\end{corollary}

 Note that, for estimation of the unknotting number,
we can use the signature and the nullity 
(see Theorem 10.1 in \cite{Mu} and Corollary 3.21 in \cite{KT})
or the sum of the absolute values of linking numbers
over all pairs of components.

 For example, the diagram $D$ in Figure \ref{fig:example}
needs at least 
$\lceil iu_{S,T}(D)/2 \rceil = \lceil 4/2 \rceil = 2$ 
RII and RIII moves 
to be deformed 
so that it has no crossing between subarcs of the same component.

 In general, for a link diagram $D$, 
the sum of the signs of all the crossings
is called the {\it writhe} and denoted by $w(D)$.
 It does not change under an RII or RIII move
but increases (resp. decreases) by one 
under an RI move creating a positive (resp. negative) crossing.
 Set
${}_{\epsilon,\delta}iu_{S,T} (D) 
= iu_{S,T}(D) + \epsilon (\dfrac{1}{2}c(D) + \delta \dfrac{3}{2}w(D))$
for a link diagram $D$,
where $\epsilon = \pm 1$, $\delta = \pm 1$
and $c(D)$ denotes the number of crossings of $D$.
 Then we have the next corollary.

\begin{corollary}\label{corollary:change2}
 The link diagram invariant 
${}_{\epsilon, +1}iu_{S,T} (D)$ (resp. ${}_{\epsilon, -1}iu_{S,T} (D)$)
increases by $2\epsilon$ 
under an RI move creating a positive (resp. negative) crossing,
decreases by $\epsilon$ 
under an RI move creating a negative (resp. positive) crossing.
 ${}_{\epsilon, \delta}iu_{S,T} (D)$ changes 
at most by $2$ under an RII or RIII move.

 Let $D_1$ and $D_2$ be link diagrams of the same oriented link.
 We need at least 
\newline
$\lceil 
|{}_{\epsilon, \delta}iu_{S,T}(D_1) - {}_{\epsilon, \delta}iu_{S,T}(D_2) | /2 
\rceil$ 
Reidemeister moves to deform $D_1$ to $D_2$.
 In particular,
when $D_2$ is a link diagram
with ${\mathcal C}(D_2) \cup \check{{\mathcal C}}(D_2) = \emptyset$,
we need at least $\lceil |{}_{\epsilon, \delta}iu_{S,T}(D_1)|/2 \rceil$ 
Reidemeister moves.
\end{corollary}

 For a knot diagram $D$,
we consider only $iu_{S,T}(D)$, ${}_{\epsilon,\delta}iu_{S,T}(D)$ with $T=(+1)$,
and we denote them by $iu_S (D)$, ${}_{\epsilon,\delta}iu_S (D)$ for short.

\begin{theorem}\label{theorem:main}
 For the diagram $U$ of the unknot in Figure \ref{fig:R8_13_1},
which can be unknotted by the sequence of $7$ Reidemeister moves in the figure,
$g(I_{lk}(U))=4$, 
max${}_{\epsilon=\pm 1, \delta=\pm 1} 
\lceil |{}_{\epsilon,\delta}iu_{(+1)}(U)/2| \rceil = 6$ and 
$\lceil |{}_{+1,+1}iu_{(-1)}(U)/2| \rceil = 7$,
where $g(I_{lk})$ denotes Hass and Nowik's invariant 
(see Section \ref{section:example} for the definition).
\end{theorem}

\begin{theorem}\label{theorem:RIII}
 There is a knot diagram which admits an RIII move
under which $iu_{(+1)}$ does not change and $iu_{(-1)}$ changes.
\end{theorem}

 Because $\epsilon(\dfrac{1}{2}c(D)+\delta\dfrac{3}{2}w(D))$ 
is unchanged under an RIII move,
the above theorem holds also 
for ${}_{\epsilon,\delta}iu_{(+1)}$ and ${}_{\epsilon,\delta}iu_{(-1)}$.
 Note that the Hass and Nowik's invariant $g(I_{lk})$
always changes by one under an RIII move.

 We prove Theorems \ref{theorem:main} and \ref{theorem:RIII} in Section \ref{section:example}.


\section{Change of the link diagram invariant under a Reidemeister move}
\label{section:change}

 Let $D$ be an oriented link diagram in the $2$-sphere, 
and $p$ a crossing of $D$.
 A (regular) {\it smoothing} at $p$ is an operation
which yields a new link diagram $D'_p$ from $D$ as below. 
 We first cut the link at the two preimage points of $p$.
 Then we obtain the four endpionts.
 We paste the four short subarcs of the link near the endpoints
in the way other than the original one
so that their orientations are connected consistently.
 See Figure \ref{fig:smoothing}.
 If we paste the four short subarcs
so that their orientations are inconsistent,
we obtain another link diagram, which we denote by $\check{D'_p}$.
 We consider $\check{D'_p}$ unoriented.
 We call this operation deforming $D$ to $\check{D'_p}$
an {\it irregular smoothing} at $p$.

 The {\it trivial $n$-component link}
is the link which has $n$ components and bounds a disjoint union of $n$ disks.
 The trivial $n$-component link admits
a link diagram with no crossings,
which we call a {\it trivial diagram}.

 Let $L$ be a link with $n$ components,
and $D$ a diagram of $L$.
 We call a sequence of Reidemeister moves and crossing changes on $D$
an {\it $X$-unknotting sequence} in this paragraph
if it deforms $D$ into a (possibly non-trivial) diagram 
of the trivial $n$-component link.
 The {\it length} of an $X$-unknotting sequence
is the number of crossing changes in it.
 The minimum length among all the $X$-uknotting sequences on $D$
is called the {\it uknotting number} of $L$.
 We denote it by $u(L)$.

 Then we define the link diagram invariants 
$iu_{S,T}$, ${}_{\epsilon,\delta}iu_{S,T}$ 
as in Section \ref{section:introduction}.

\begin{figure}[htbp]
\begin{center}
\includegraphics[width=10cm]{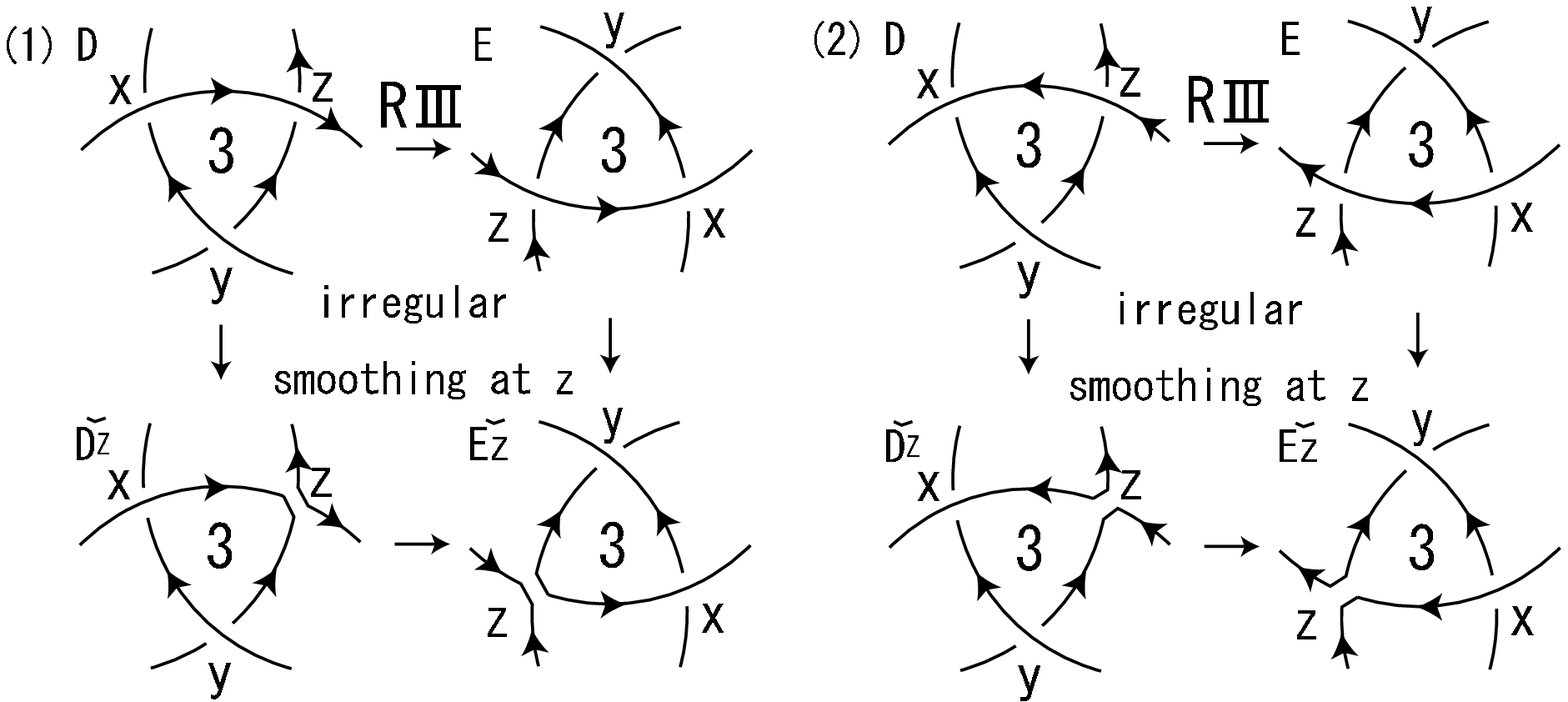}
\end{center}
\caption{}
\label{fig:IrregularTrigon}
\end{figure}

\begin{proof}[Proof of Theorem \ref{theorem:change}]
 The proof is very similar to the arguments 
in Section 2 in \cite{HN1} and in Section 2 in \cite{HH2}.
 Let $D, E$ be link diagrams of a link $L$
such that $E$ is obtained from $D$ by a Reidemeister move.

 First, 
we suppose that $E$ is obtained from $D$ by an RI move creating a crossing $a$.
 Let $L_i$ be a component of $L$ 
such that $a$ is a crossing between two subarcs of $L_i$.
 When the element $s_{ii}$ of the symmetric matrix $S$ is $+1$,
the same argument as in \cite{HH2} works, and we omit it.
 We consider the case where $s_{ii}=-1$.
 Then the link $\check{E_a}$ obtained from $E$ by an irregular smoothing at $a$
is of the same link type as the original link $L$,
and hence $u(\check{E_a}) = u(L)$.
 Then the contribution of $a$ to $iu_{S,T}(E)$ is $\pm(u(\check{E_a})-u(L))= 0$.
 The contribution of any other crossing $x$ to $iu_{S,T}$ is unchanged
since the RI move shows 
that $D_x$ and $E_x$ are the same link.
 Thus the RI move does not change $iu_{S,T}$, i.e., $iu_{S,T}(D)=iu_{S,T}(E)$.
 When $s_{ii}=0$, 
the crossing $a$ does not contributes to $iu_{S,T}(E)$,
and $iu_{S,T}$ does not change.

 We consider the case
where $E$ is obtained from $D$ by an RII move creating a bigon face.
 Let $x$ and $y$ be 
the positive and negative crossings at the corners of the bigon,
and $L_i, L_j$ components of $L$
such that the edges of the bigon are a subarc of $L_i$ and a subarc of $L_j$.
 It is possible that $i=j$.
 When $s_{ij}=+1$, the proof is the same as in Section 2 in \cite{HH2},
and we omit it.
 When $s_{ij}=0$, it is clear that $iu_{S,T}$ does not change.
 Hence we assume that $s_{ij}=-1$.
 If the RII move is matched,
then $\check{E_x}$ and $\check{E_y}$ are the same link.
 Hence $u(\check{E_x})=u(\check{E_y})$ and 
$|iu_{S,T}(E)-iu_{S,T}(D)| 
= |t_{ij}(|u(\check{E_x})-u(L)| - |u(\check{E_y})-u(L)|)| = 0$.
 If the RII move is unmatched,
then a crossing change on $\check{E_x}$ yields $\check{E_y}$, 
and hence their unknotting numbers differ by at most one,
i.e., $|u(\check{E_x})-u(\check{E_y})| \le 1$. 
 Hence 
$|iu_{S,T}(E)-iu_{S,T}(D)| 
= ||u(\check{E_x})-u(L)| - |u(\check{E_y})-u(L)||
\le |(u(\check{E_x})-u(L))-(u(\check{E_y})-u(L))|
= |u(\check{E_x})-u(\check{E_y})|
\le 1$.

 We consider the case where $E$ is obtained from $D$ by an RIII move.
 See Figure \ref{fig:IrregularTrigon} for typical examples.
 In $D$ and $E$, let $x$ be the crossing 
between the top and the middle strands of the trigonal face
where the RIII move is applied.
 Let $L_i$ and $L_j$ be the components of $L$
such that $L_i$ contains the top strand and $L_j$ contains the middle strand.
 Then $D_x$ and $E_x$ are the same link when $s_{ij}=+1$. 
 Also, irregular smoothing operations on $D$ and $E$ at $x$
yields the same link $\check{D_x}$ and $\check{E_x}$
when $s_{ij}=-1$.
 Hence the contribution of $x$ to $iu_{S,T}$ is unchanged under the RIII move.
 The same is true for the crossing $y$ 
between the bottom and the middle strands.
 Let $z$ be the crossing 
between the top strand in $L_m$ and the bottom strand in $L_n$.
 The proof is the same as in Section 2 in \cite{HH2}.
 We recall it in the case where $s_{mn}=-1$.
 $\check{D_z}$ and $\check{E_z}$ differ 
by Reidemeister moves and at most two crossing changes,
and hence $|u(\check{E_z})-u(\check{D_z})| \le 2$. 
 Thus
$|iu_{S,T}(E)-iu_{S,T}(D)|
= ||u(\check{E_z})-u(L)| - |u(\check{D_z})-u(L)||
\le |(u(\check{E_z})-u(L))-(u(\check{D_z})-u(L))|
= |u(\check{E_z})-u(\check{D_z})|
\le 2$.
\end{proof}


\section{Examples}\label{section:example}

 In this section, we prove Theorems \ref{theorem:main} and \ref{theorem:RIII}.

 We first recall the definition of Hass and Nowik's invariant $g(I_{lk})$
defined in \cite{HN1} and \cite{HN2}.
 Let ${\mathbb G}_{\mathbb Z}$ be 
the free abelian group with basis $\{ X_n , Y_n \}_{n \in {\mathbb Z}}$.
 The invariant $I_{lk}$ 
assigns an element of ${\mathbb G}_{\mathbb Z}$ to a knot diagram.
 Let $D$ be an oriented knot diagram.
 Then $I_{lk}$ is defined by the formula
$I_{lk}(D)=
\sum_{p \in {\mathcal C}^+ (D)} X_{lk(D_p)}
+
\sum_{m \in {\mathcal C}^- (D)} Y_{lk(D_m)}$,
where ${\mathcal C}^+ (D)$ (resp. ${\mathcal C}^- (D)$) is 
the set of all the positive (resp. negative) crossings of $D$
and $lk$ is the linking number.
 Let $g : {\mathbb G}_{\mathbb Z} \rightarrow {\mathbb Z}$ 
be the homomorphism 
defined by $g(X_n)=|n|+1$ and $g(Y_n)=-|n|-1$.
 Then the numerical invariant $g(I_{lk}(D))$ of a knot diagram $D$
changes at most by one under a Reidemeister move.
 This invariant can be decomposed as
$g(I_{lk}(D))=g_0(I_{lk}(D))+w(D)$,
where $w(D)$ is the writhe
and the homomorphism $g_0 : {\mathbb G}_{\mathbb Z} \rightarrow {\mathbb Z}$ 
is defined by $g_0(X_n)=|n|$ and $g_0(Y_n)=-|n|$.
 It can be easily seen 
that $g_0(I_{lk}(D))$ does not change under an RI move
and the change of $g_0(I_{lk}(D))$ under an RII or RIII move is the same
as that of $g(I_{lk}(D))$.

\begin{figure}[htbp]
\begin{center}
\includegraphics[width=12cm]{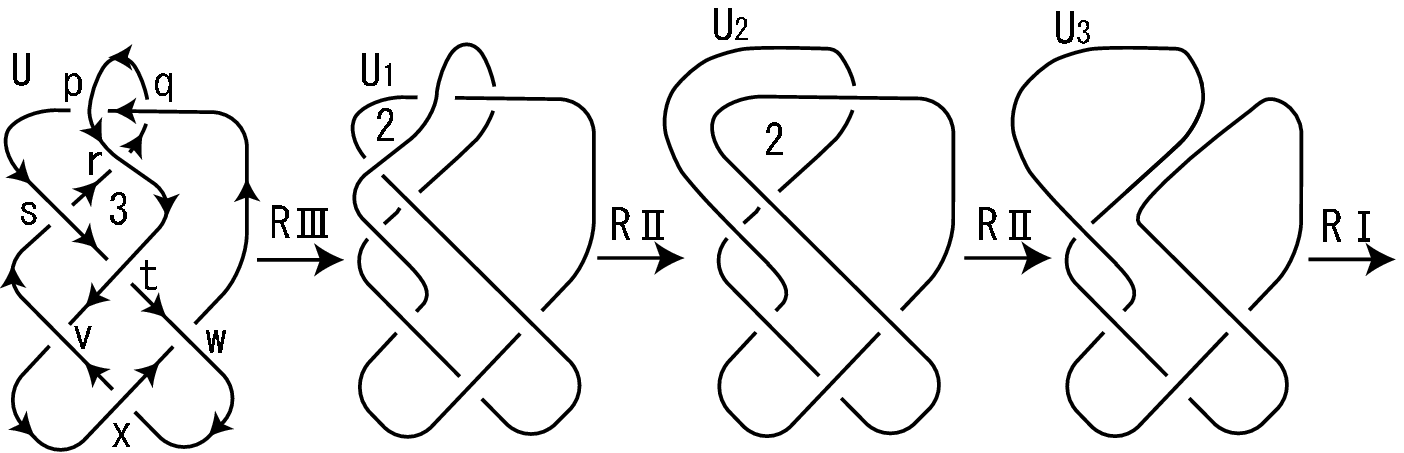}
\end{center}
\caption{}
\label{fig:R8_13_1}
\end{figure}

\begin{proof}[Proof of Theorem \ref{theorem:main}]
 We consider the diagram $U$ in Figure \ref{fig:R8_13_1},
which represents the trivial knot
unknotted by the sequence of Reidemeister moves in Figure \ref{fig:R8_13_1}.
 We can unknot the last diagram by four RI moves.

\begin{figure}[htbp]
 \begin{minipage} {0.4\hsize}
  \begin{center}
   \includegraphics[width=55mm]{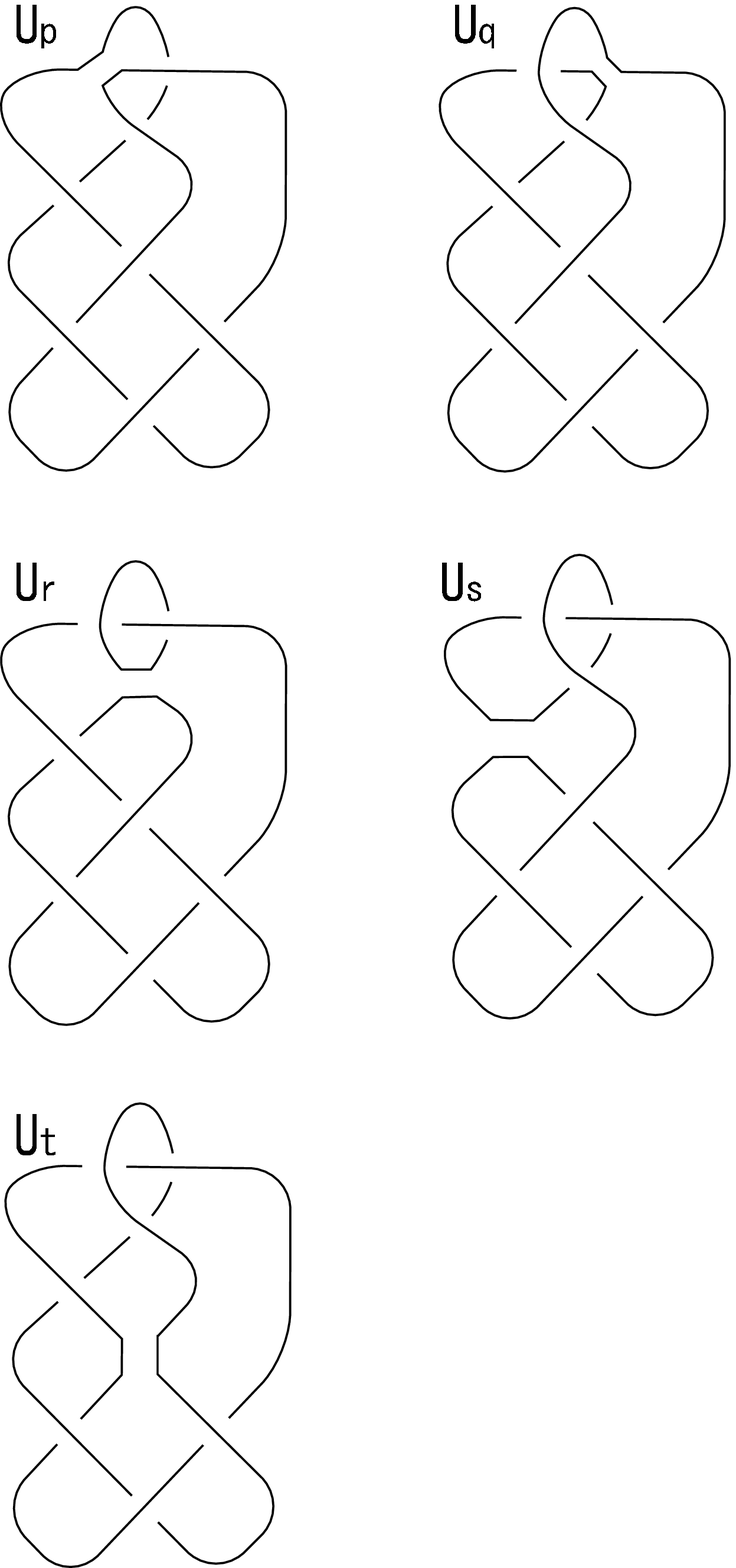}
  \end{center}
  \caption{}
  \label{fig:smooth8_13_1}
 \end{minipage}
 \begin{minipage} {0.4\hsize}
  \begin{center}
   \includegraphics[width=60mm]{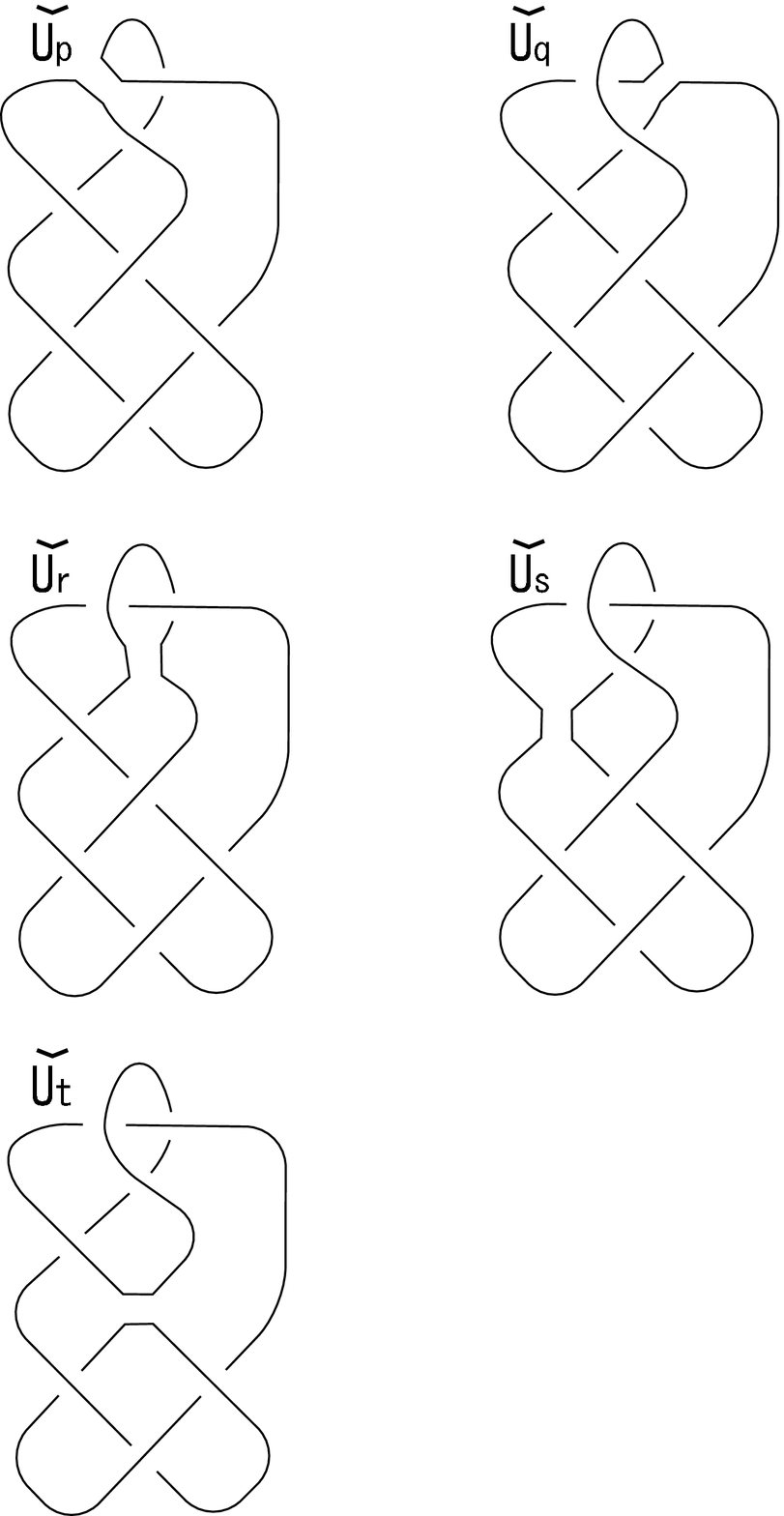}
  \end{center}
  \caption{}
  \label{fig:irregular8_13_1}
 \end{minipage} 
\end{figure}

 For this diagram,
$g_0 (I_{lk}(U))= 0$,
$iu_{(+1)}(U)/2 = 1/2$ and 
$iu_{(-1)}(U)/2 = 3/2$.

 In fact, $U$ has $8$ crossings
$p$ through $t$ and $v$ through $x$.
 We obtain the links in Figure \ref{fig:smooth8_13_1}
from $U$ by a regular smoothing operation.
 Each of $U_p$, $U_q$ and $U_s$ is the $(2,4)$-torus link 
with one component furnished the reversed orientation.
 Hence $|lk(U_p)|=|lk(U_q)|=|lk(U_s)|=2$
and $u(U_p)=u(U_q)=u(U_s)=2$.
 $U_r$ is a composite of the twist knot $5_2$ 
in the Rolfsen table \cite{Ro}
and the Hopf link.
 Hence $|lk(U_r)|=1$ and $u(U_r)=2$.
 Note that 
we need a crossing change at a crossing between distinct components
for unlinking,
and another crossing change 
at a crossing between subarcs of the twist knot component
for unknot the twist knot. 
 $U_t$ is the Hopf link,
and hence $|lk(U_t)|=1$ and $u(U_t)=1$.
 $U_v$, $U_w$ and $U_x$ are the trivial $2$-component links, and $|lk|=u=0$ for these links.
 Since ${\rm sign}\,p={\rm sign}\,q=-1$, and $+1$ for the other crossings, 
$g_0(I_{lk}(U)) = -2-2+1+2+1+0+0+0 = 0$,
and 
$iu_{(+1)}(U) = -2-2+2+2+1+0+0+0 = 1$.

 We obtain the knots in Figure \ref{fig:irregular8_13_1}
from $U$ by an irregular smoothing operation.
 $\check{U_p}$, $\check{U_q}$ and $\check{U_t}$ are the trefoil knots,
and $u(\check{U_p})=u(\check{U_q})=u(\check{U_t})=1$.
 $\check{U_r}$ is the $7_4$ knot, and $u(\check{U_r})=2$ by \cite{L}.
 $\check{U_s}$ is the star knot $5_2$, and $u(\check{U_s})=2$ (see \cite{Mu}).
 $\check{U_v}$, $\check{U_w}$ and $\check{U_x}$ are the trivial knots,
and $u(\check{U_v})=u(\check{U_w})=u(\check{U_x})=0$.
 Hence $iu_{(-1)}(U)=-1-1+2+2+1+0+0+0=3$.

 Since $c(U)=8$ and $w(D)=4$, 
the expression $\epsilon (\dfrac{1}{2}c(D) + \delta \dfrac{3}{2}w(D))$ 
takes the maximum value $10$
when $\epsilon=+1$ and $\delta=+1$.
 Because $iu_{(+1)}(U)$ and $iu_{(-1)}(U)$ are positive,
$|{}_{\epsilon,\delta}iu_{(+1)}(U)|$
and 
$|{}_{\epsilon,\delta}iu_{(-1)}(U)|$
are the largest when $\epsilon=+1$ and $\delta=+1$.
 Then easy calculations show the theorem.
\end{proof}

\begin{remark}
 The length of the unknotting sequence of Reidemeister moves 
in Figure \ref{fig:R8_13_1} is $7$.
 Without using ${}_{+1, +1}iu_{(-1)}(U)$,
we can easily see that this is minimal as below.
 Since the writhe $w(U)=4$, we need four RI moves for unknotting.
 These RI moves decrease the number of crossings of $U$ by four.
 Because $U$ has eight crossings,
we need to delete four more crossings.
 Hence at least two more RI and RII moves are necessary for unknotting.
 We need at least one more Reidemeister move
since $U$ has neither a monogon face nor a bigon face
on which we cay apply an RI or RII move to decrease the number of crossings.

 $g(I_{lk})$ does not change 
under the second RII move in Figure \ref{fig:R8_13_1}
because the move is unmatched.
 Under the sequence of Reidemeister moves in Figure \ref{fig:R8_13_1},
$iu_{(+1)}$ varies 
$1 \rightarrow -1 \rightarrow 0 \rightarrow 0 
\rightarrow 0 \rightarrow 0 \rightarrow 0 \rightarrow 0$,
$iu_{(-1)}$ varies 
$3 \rightarrow 1 \rightarrow 1 \rightarrow 0 
\rightarrow 0 \rightarrow 0 \rightarrow 0 \rightarrow 0$,
${}_{+1,+1}iu_{(+1)}$ varies 
$11 \rightarrow 9 \rightarrow 9 \rightarrow 8 
\rightarrow 6 \rightarrow 4 \rightarrow 2 \rightarrow 0$,
${}_{+1,+1}iu_{(-1)}$ varies 
$13 \rightarrow 11 \rightarrow 10 \rightarrow 8 
\rightarrow 6 \rightarrow 4 \rightarrow 2 \rightarrow 0$.
 So, ${}_{+1,+1}iu_{(-1)}$ 
is the most sensitive to each Reidemeister move
in this case.
\end{remark}

\begin{figure}[htbp]
 \begin{minipage} {0.4\hsize}
  \begin{center}
   \includegraphics[width=70mm]{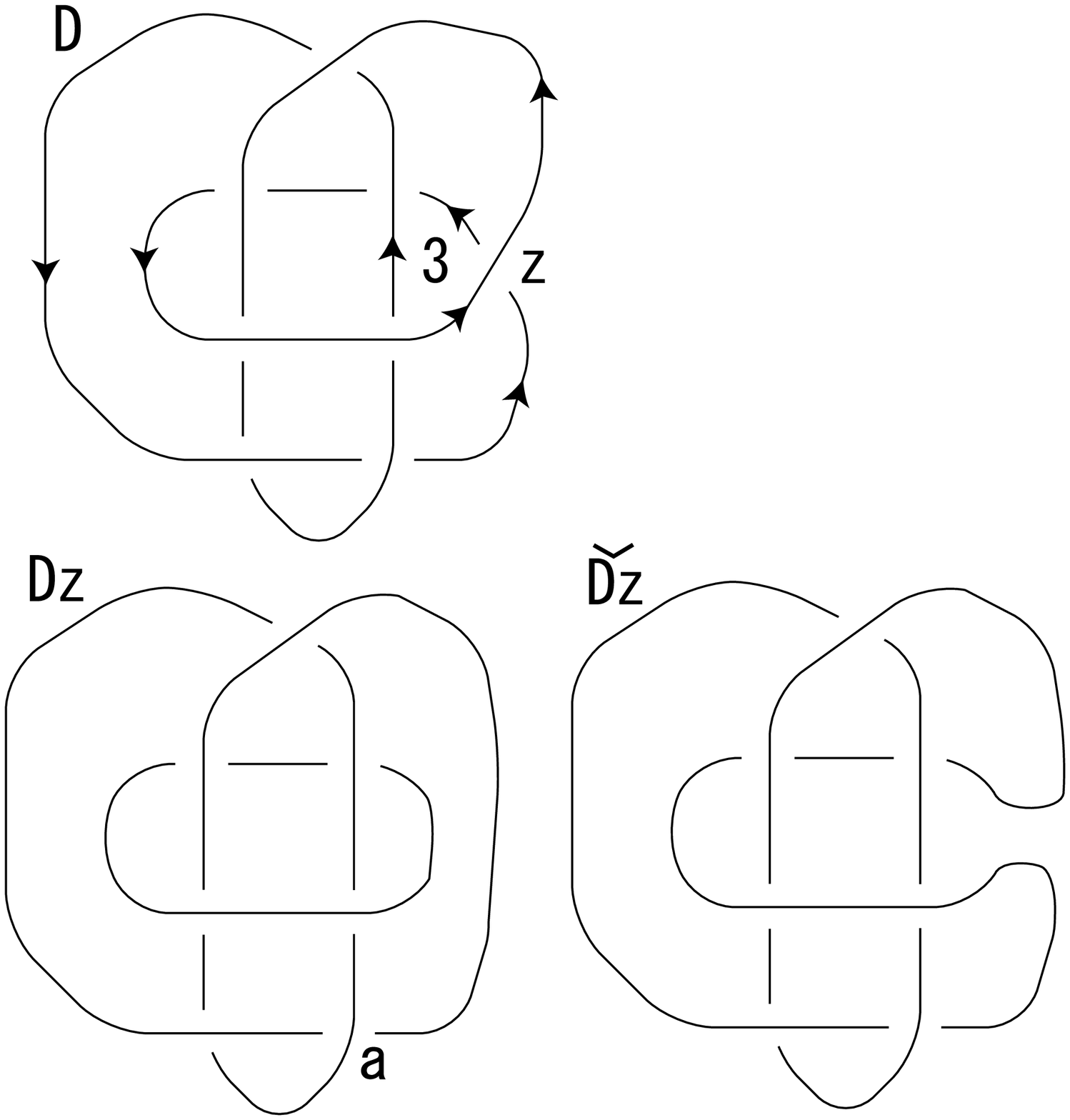}
  \end{center}
  \caption{}
  \label{fig:smooth8_17_13}
 \end{minipage}
 \begin{minipage} {0.4\hsize}
  \begin{center}
   \includegraphics[width=70mm]{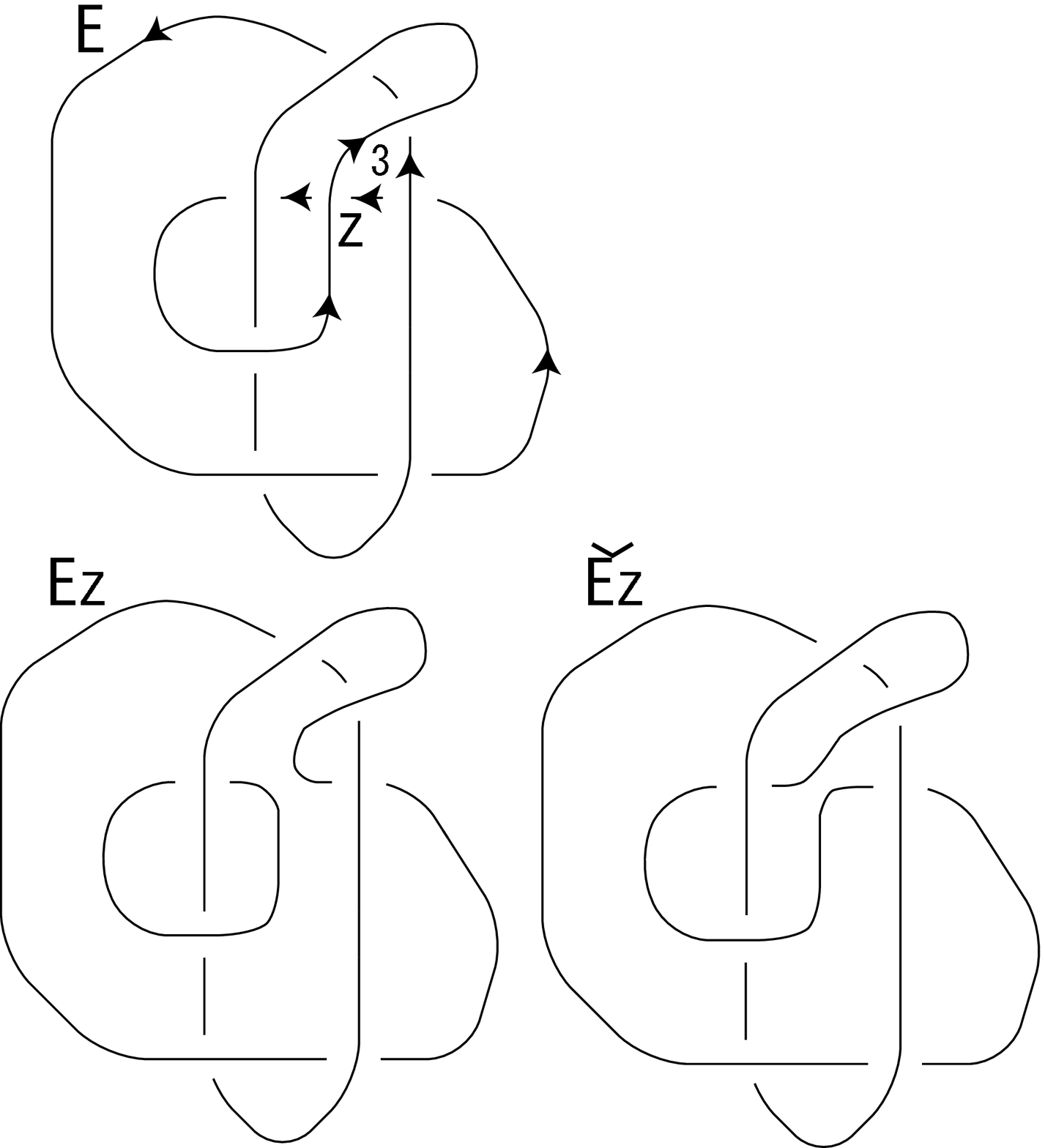}
  \end{center}
  \caption{}
  \label{fig:R8_17_13}
 \end{minipage} 
\end{figure}

\begin{figure}[htbp]
 \begin{minipage} {0.4\hsize}
  \begin{center}
   \includegraphics[width=60mm]{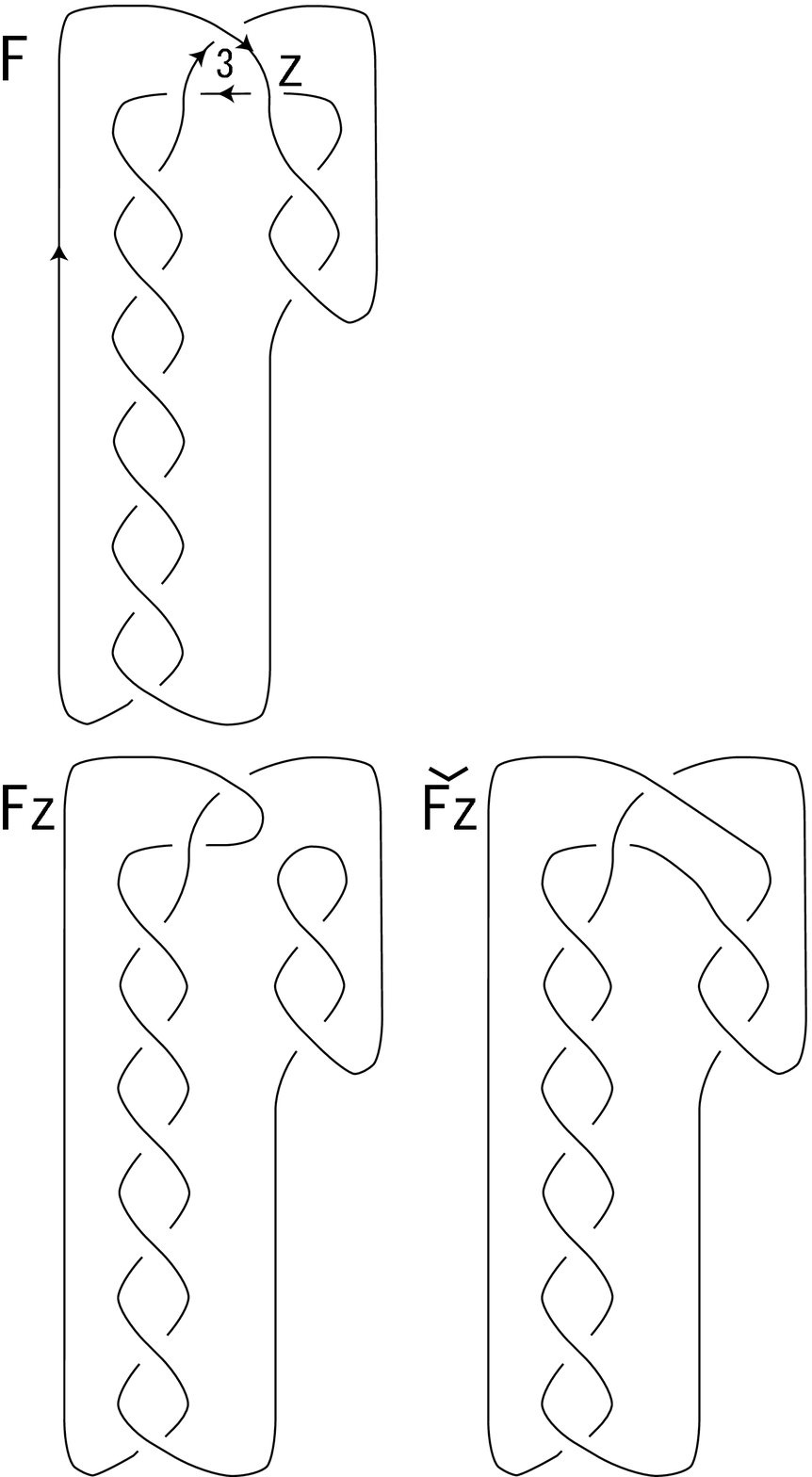}
  \end{center}
  \caption{}
  \label{fig:twistknotsmooth}
 \end{minipage}
 \begin{minipage} {0.4\hsize}
  \begin{center}
   \includegraphics[width=60mm]{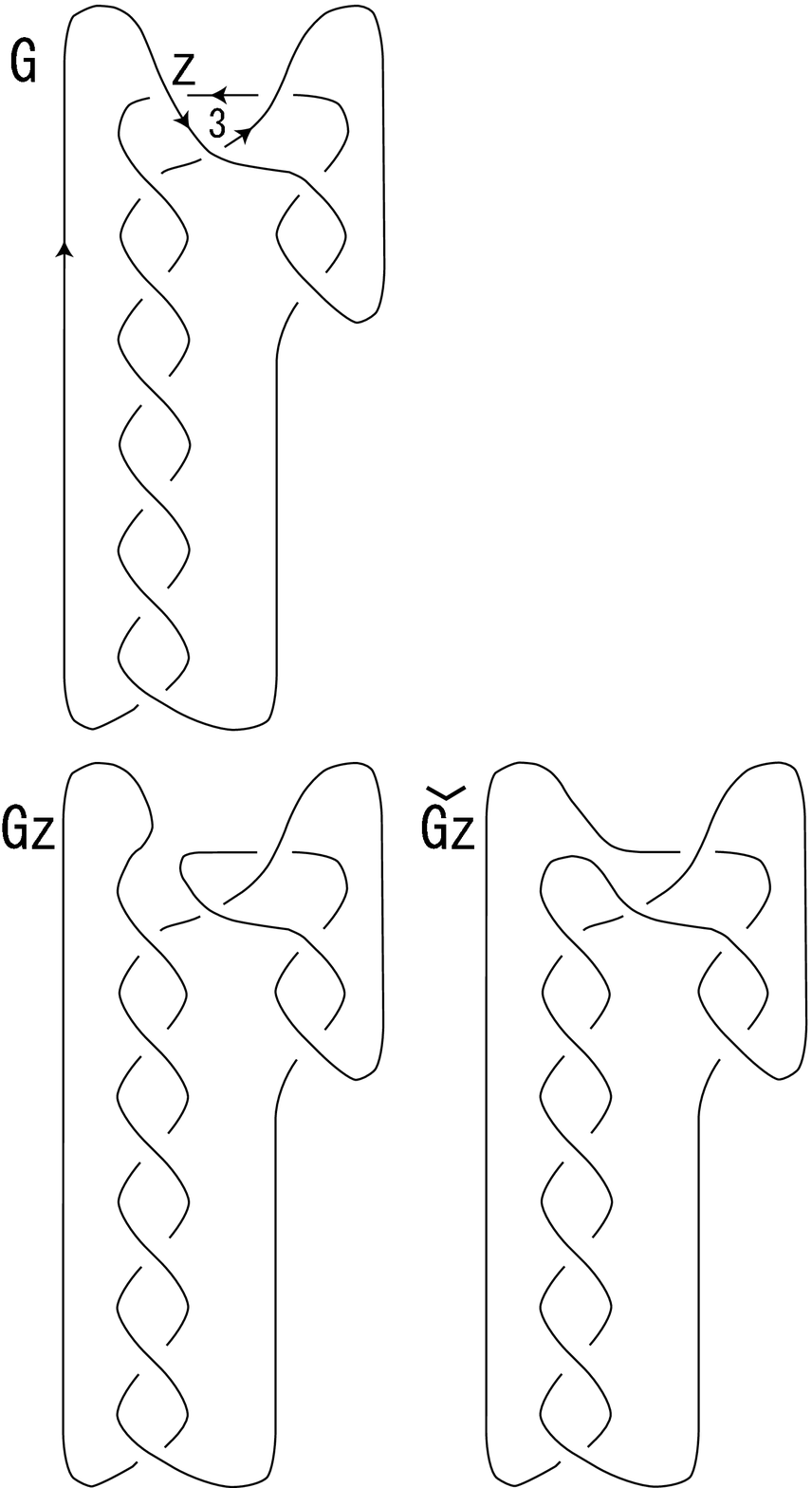}
  \end{center}
  \caption{}
  \label{fig:Rtwistknot}
 \end{minipage} 
\end{figure}

\begin{proof}[Proof of Theorem \ref{theorem:RIII}]
 We give two examples:
the first one is an RIII move on a diagram of the trivial knot
in Figures \ref{fig:smooth8_17_13} and \ref{fig:R8_17_13}
under which $iu_{(+1)}$ does not change
and $iu_{(-1)}$ changes by one,
and the second one is an RIII move on a diagram of a twist knot
in Figures \ref{fig:twistknotsmooth} and \ref{fig:Rtwistknot}
under which $iu_{(+1)}$ does not change
and $iu_{(-1)}$ changes by two.
 The RIII moves are performed at the trigonal faces labeled $3$.
 In these diagrams, let $z$ be the crossing 
between the top and the bottom strands of the trigonal face.

 The diagram $D$ in Figure \ref{fig:smooth8_17_13} is deformed
into the diagram $E$ in Figure \ref{fig:R8_17_13}.
 As shown in the argument in the proof of Theorem \ref{theorem:change}
in Section \ref{section:change},
the changes of $iu_{(+1)}$ and $iu_{(-1)}$ are calculated as below.
\newline
$|iu_{(+1)}(E)-iu_{(+1)}(D)|
= ||u(E_z)-u(L)| - |u(D_z)-u(L)||
= ||1-0|-|1-0||=0$
\newline
$|iu_{(-1)}(E)-iu_{(-1)}(D)|
= ||u(\check{E_z})-u(L)| - |u(\check{D_z})-u(L)||
= ||1-0|-|0-0||=1$
\newline
 Note that $u(D_z)=1$ 
since a crossing change at the crossing $a$ in Figure \ref{fig:R8_17_13}
yields the trivial $2$-component link.
 $u(E_z)=1$, $u(\check{D_z})=1$ and $u(\check{E_z})=0$
because $E_z$ is the Hopf link, $\check{D_z}$ the twist knot $5_2$
and $\check{E_z}$ the trivial knot.

 The diagram $F$ in Figure \ref{fig:twistknotsmooth} is deformed
into the diagram $G$ in Figure \ref{fig:Rtwistknot}.
 The changes of $iu_{(+1)}$ and $iu_{(-1)}$ are calculated as below.
\newline
$|iu_{(+1)}(G)-iu_{(+1)}(F)|
= ||u(G_z)-u(L)| - |u(F_z)-u(L)||
= ||4-1|-|4-1||=0$
\newline
$|iu_{(-1)}(G)-iu_{(-1)}(F)|
= ||u(\check{G_z})-u(L)| - |u(\check{F_z})-u(L)||
= ||1-1|-|3-1||=2$
\newline
 Note that $u(\check{F_z})=3$ 
because $\check{F_z}$ is the knot $10_2$ in the Rolfsen table \cite{Ro}
and its signature is $-6$ (see \cite{Mu}).
 $u(G_z)=4$
since $G_z$ is the composite of the $(2,6)$-torus link and the trefoil knot.
 (We need at least $3$ crossing change 
at crossings between distinct components
to split the two components
because $lk(G_z)=3$,
and at least one crossing change between subarcs of the trefoil knot component
to unknot this component.)
 $u(F_z)=4$ and $u(\check{G_z})=1$
because $F_z$ is the $(2,8)$-torus link
and $\check{G_z}$ the trefoil knot.
\end{proof}



\bibliographystyle{amsplain}

\medskip

\noindent
Chuichiro Hayashi: 
Department of Mathematical and Physical Sciences,
Faculty of Science, Japan Women's University,
2-8-1 Mejirodai, Bunkyo-ku, Tokyo, 112-8681, Japan.
hayashic@fc.jwu.ac.jp

\vspace{3mm}
\noindent
Miwa Hayashi:
Department of Mathematical and Physical Sciences,
Faculty of Science, Japan Women's University,
2-8-1 Mejirodai, Bunkyo-ku, Tokyo, 112-8681, Japan.
\newline
miwakura@fc.jwu.ac.jp

\end{document}